\documentstyle[11pt]{article}
\renewcommand{\baselinestretch}{1.2}
\begin{document}

\newtheorem{theorem}[section]{Theorem}
\newtheorem{lemma}[section]{Lemma}
\newtheorem{corollary}[section]{Corollary}
\newtheorem{proposition}[section]{Proposition}
\newenvironment{proof}{{\bf Proof.}}{\hfill{ }\fbox{$^{\mbox{ }}$}}
\renewcommand{\baselinestretch}{1.5}\newtheorem{conjecture}[section]{Conjecture}

\begin{center}
{\Large {\bf On Arithmetic Progressions of Cycle Lengths in Graphs}}
\end{center}

\begin{center}
Jacques Verstra\"{e}te \\
Department of Pure Mathematics and Mathematical Statistics \\
Centre for Mathematical Sciences \\
Wilberforce Road, Cambridge CB3 OWB \\
England August 1999. \\
{\sf jbav@microsoft.com}
\end{center}

\bigskip
\bigskip

\begin{center}
{\bf Abstract}

\medskip

\renewcommand{\baselinestretch}{1}
\parbox{4in}{\small 
A recently posed question of H\"{a}ggkvist and Scott's asked whether or not there exists a constant $c$ such that if $G$ is a graph of minimum degree $ck$ then $G$ contains cycles of $k$ consecutive even lengths. In this paper we answer the question by proving that for $k \geq 2$, a bipartite graph of average degree at least $4k$ and girth $g$ contains cycles of $(g/2-1)k$ consecutive even lengths. We also obtain a short proof of the theorem of Bondy and Simonovits, that a graph of order $n$ and size at least $8(k-1)n^{1 + 1/k}$ has a cycle of length $2k$.}
\end{center}

\bigskip

Erd\H{o}s and Burr~\cite{Er1} conjectured that for every odd number $k$, there is a constant $c_{k}$ such that for every natural number $m$, every graph of average degree at least $c_{k}$ contains a cycle of length $m$ modulo $k$.  Erd\H{o}s and Burr~\cite{Er1} settled their conjecture in the case $m = 2$ and Robertson (see~\cite{Er1}) settled the case $m = 0$.
The full conjecture was resolved by Bollob\'{a}s~\cite{BB}, who proved the conjecture with $c_{k} = 2[(k + 1)^{k} - 1]/k$. In this paper, we show that 
$c_{k} = 8k$ will do. 
 Thomassen~\cite{Th2} later showed cycles of all even lengths modulo $k$ are obtained under the hypothesis that the average degree is at least $4k(k+1)$, without requiring $k$ to be odd.
Thomassen~\cite{Th1} also proved that if $G$ is a graph of minimum degree at least three and girth at least $2(k^{2} + 1)(3\cdot 2^{k^{2}+1} + (k^{2} + 1)^{2} - 1)$, then $G$ contains cycles of all even lengths modulo $k$. 

\medskip

Bondy and Vince~\cite{BV} proved that in a graph in which all but at most two vertices have degree at least three, there exist two cycles whose lengths differ by at most two. This answered a conjecture of Erd\H{o}s and was also studied by H\"{a}ggkvist and Scott~\cite{HS2}.
Recently, H\"{a}ggkvist and Scott~\cite{HS1}, considered extending this to considering arithmetic progressions of cycle lengths in graphs. H\"{a}ggkvist and Scott~\cite{HS1} proved that if $G$ is a graph 
of minimum degree at least $300k^{2}$ then $G$ contains $k + 1$ consecutive 
even cycle lengths. The same authors asked if a linear bound on the 
minimum degree is possible. In this paper, we answer the question of H\"{a}ggkvist and Scott in the following theorem.
\begin{theorem}\label{cyc}
 Let $k \geq 2$ be a natural number and $G$ a bipartite graph of average degree at least $4k$ and girth $g$. Then there exist cycles of $(g/2 - 1)k$ consecutive even lengths in $G$. Moreover, the shortest of these cycles has length at most twice the radius of $G$. 
\end{theorem}
This result generalises the above-mentioned result of Bollob\'{a}s and partly generalises that of Thomassen, insofar as Thomassen's result is valid for graphs of minimum degree at least three, whereas the result above requires average degree at least eight. The graph $K_{k,n-k}$ shows that we require the average degree to be at least about $2k$ to ensure the conclusion of Theorem 1.

\medskip
 
The following lemma lies at the heart of the proof of Theorem 1. It was originally inspired by methods used by Gy\'{a}rf\'{a}s, Koml\'{o}s and Szemer\'{e}di~\cite{GKS}. Whilst this paper was being written, the lemma was discovered to be implicit in a lemma of Bondy and Simonovits~\cite{BS}. Nevertheless, the proof is short and is retained here for completeness. 

\begin{lemma} Let $H$ be a graph comprising a cycle with a chord. Let
$(A,B)$ be a non-trivial partition of $V(H)$. Then $H$ contains $A$--$B$
paths of every length less than $|H|$, unless $H$ is bipartite with bipartition
$(A,B)$.
\end{lemma}

 \begin{proof} Label the vertices of the cycle $0,1,\ldots,n-1$
where $n=|H|$. Suppose $H$ does not contain $A$--$B$ paths of every length
less than $n$, and let $m$ be the smallest integer for which there is no
$A$--$B$ path of length $m$ not using the chord; $m>1$ since $(A,B)$ is a
non-trivial partition of $V(H)$. We remark also that $m \leq n/2$, or $H$ would contain $A$--$B$ paths of all lengths less than $n$.

\smallskip

Now $\chi(j)=\chi(j+m)$ for every $j\in V(H)$, where $\chi$ is the
characteristic function of $A$ (label arithmetic is modulo $n$). Let $d=\hbox{hcf}(n,m)$. Then there are integers $p$ and $q$ such that
$pm+qn=d$; hence $\chi(j)=\chi(j+d)$ for every $j$. But then there is no
$A$--$B$ path of length $d$ round the cycle; thus $d=m$ and $m\mid n$. In
particular, $A$--$B$ paths of every length less than $m$ exist by the
definition of $m$, so $A$--$B$ paths of every length other than
multiples of $m$ exist by periodicity of $\chi$.

\smallskip

We find paths of the remaining lengths $km$, $1\le k \le n/m -1$, using
the chord. Suppose first that the chord joins two vertices within distance $m$
on the cycle, say $0$ and $r$ where $1<r\le m$. There exist $A$--$B$ paths of
length $m+r-1$ round the cycle; thus $\chi(j)\ne\chi(j+m+r-1)$ for some $j$,
$-m<j\le0$. But $j+m+r-1\ge r$, so the path
$j,j+1,\ldots,0,r,r+1,\ldots,j+km+r-1$ is an $A$--$B$ path of length $km$
provided $j+km+r-1 < n+j$, which holds for all the desired
$k\le n/m -1$.

\smallskip

So we may suppose the chord is $0r$, where $m< r < n-m$. Let $-m<j<0$ and
consider the paths $j,j+1,\ldots,0,r,r-1,\ldots,r-j-m+1$ and
$m+j,m+j-1,\ldots,0,r,r+1,\ldots,r-j-1$, of length $m$. If either of them is an
$A$--$B$ path we can extend it, by $m$ vertices at a time, to $A$--$B$ paths of
lengths $km$, $k\ge 1$, until the number of unused vertices in the two arcs
defined by the chord is less than $m$ in each arc. At this point $km+1\ge
n-2(m-1)$, and as $m\mid n$, $km=n-m$ as desired. Likewise, if
either of the two paths $0,r,r-1,\ldots,r-m+1$ and $0,r,r+1,r+m-1$ is an
$A$--$B$ path then $H$ contains paths of all lengths less than $|H|$.

\smallskip

Thus it follows, as $\chi(j)=\chi(m+j)$, that for $-m<j<0$ we have
$\chi(r-j-m+1)= \chi( r-j-1)$ and that $\chi(r-m+1)=\chi(r+m-1)$, implying
$\chi(r-j+1)= \chi( r-j-1)$ and $\chi(r+m+1)=\chi(r+m-1)$. So
$\chi(v+2)=\chi(v)$ for all $r\le v< r+m$, and so for all $v\in
V(H)$. Hence $m=2$.

\smallskip

We conclude, therefore, that $|H|$ is even and the vertices of the cycle are
alternately in $A$ and in $B$. It is immediately seen that, under these
circumstances, if the chord joins two vertices in the same class then $H$
contains $A$--$B$ paths of all lengths less than $|H|$. Consequently, the chord joins
$A$ to $B$, so $H$ is bipartite, with bipartition $(A,B)$. 
\end{proof} \medskip

\begin{lemma} Let $k \geq 2$ be a natural number and let $G$ be a graph of average degree at least $2k$ and girth $g$. Then $G$ contains a cycle of length at least $(g-2) k  + 2$, with at least one chord. 
\end{lemma}

\begin{proof} It is easily seen that a graph $G$ of average degree at least $2k$ contains a subgraph $H$ of minimum degree at least $k + 1$. If $P$ is a longest path in $H$, then an endvertex $v$ of $P$ has all its neighbours on $P$. Some neighbour $u$ of $v$ is at distance at least $(g-2)k + 1$ from $v$ on $P$. Hence $P + uv$ is a cycle of length at least $(g-2) k + 2$. As $k + 1 \geq 3$, this cycle has at least one chord.
\end{proof} 

\medskip \medskip

\noindent{\bf Proof of Theorem~\ref{cyc}.} We may assume that $G$ is connected and let the radius of $G$ be $\mbox{rad}(G)$. Choose a central vertex $v_{0} \in V(G)$, and let $V_{i}$ denote the set of vertices a distance $i$ from $v_{0}$ in $G$. Then there exists $l$ such that $V_{l} \cup V_{l+1}$ 
spans a graph $G'$ with at least $k|V_{l} \cup V_{l+1}|$ edges. By Lemma 3, find $H \subset G'$ comprising a cycle, of length at least $(g - 2)k + 2$, with a chord. Let $T'$ be a minimal subtree of $T$, restricted to $\bigcup_{i \leq l} V_{i}$, such that $T'$ contains $V(H) \cap V_{l}$. The minimality of $T'$ ensures that it branches at its root. Now let $A$ be the set of vertices of $H$ in one of these branches and let $B = V(H) \backslash A$. By Lemma 2, and as $(A,B)$ is not the bipartition of $H$, there are $A$--$B$ paths of all lengths up to $(g-2)k + 1$, all disjoint from $T' - \mbox{end}(T')$. Each $A$--$B$ path of even length $s$, together with a subpath of $T'$ between the ends of such a path, gives rise to a cycle of length $s + 2r$, where $r$ is the distance from $V_{l}$ to the root of $T'$. Note that, as $G$ is bipartite, all paths of even length with one end in $A$ have their other end in $V_{l}$. This gives cycles $C_{2r+2}, C_{2r + 4}, .., C_{2r + (g-2)k}$, of $(g/2 - 1)k$ consecutive even lengths, and since $v_{0}$ is a central vertex, $2r + 2 \leq \mbox{rad}(G)$, as required. 
\endproof

\bigskip

We define the even girth of a graph $G$ to be the length of a shortest even cycle in $G$. Theorem~\ref{cyc} easily extends to general graphs, as is shown by the following corollary.

\begin{corollary} Let $k \geq 2$ be a natural number, and let $G$ be a graph of average degree at least $8k$ and even girth $g$. Then there are cycles of $(g/2 - 1)k$ consecutive even lengths in $G$. 
\end{corollary}

\begin{proof} This follows from the observation that a graph of average degree at least $8k$ has a spanning bipartite subgraph of average degree at least $4k$, and then applying Theorem 1 to this bipartite subgraph. 
\end{proof} \medskip

In the case of graphs of average degree at least $6k$, we may also argue as follows. Given a vertex $v_{0}$, let $V_{l}$ denote the vertices a distance $l$ from $v_{0}$. Then either there exists $l$ such that $V_{l} \cup V_{l+1}$ spans a bipartite graph of average degree at least $2k$, or there exists $l$ such that $V_{l}$ spans a graph of average degree at least $2k$. In the former case, the method of Theorem 1 gives cycles of all even lengths in an integer interval of form $[2r+1,2r+(g-2)k+1]$ and the latter case gives (also following the proof of Theorem 1) cycles of all odd lengths in an integer interval of form $[2r+1,2r+(g-2)k+1]$. So we have the following theorem:

\begin{theorem} Let $k \geq 2$ be a natural number, and let $G$ be a graph of average degree at least $6k$ and girth $g$. Then, for some odd number $r \geq 3$, there exist cycles of all even lengths or all odd lengths in the interval $[r,r+(g-2)k]$.
\end{theorem}

The above result recalls the result of Bondy and Vince~\cite{BV}, 
that if $G$ is a graph with at most two vertices of degree at most two, then $G$ contains cycles of two consecutive lengths or two consecutive even lengths. 
 This was proved using a technique of Thomassen and Toft~\cite{TT}.
 In comparison, note that Theorem 5 requires average degree at least $6k$ where $k \geq 2$. Therefore, to ensure two cycles of consecutive lengths or consecutive even lengths, we require average degree at least twelve, which is higher than what is required in the context of Bondy and Vince's results. 

\medskip

The following result is proved in the same was as Theorem 1:

\begin{corollary}
Let $k \geq 2$ be a natural number, and suppose that $G$ has chromatic 
number at least $2k + 2$ and girth at least $g$. 
Then $G$ contains cycles of $k(g  - 2)$ consecutive 
lengths.
\end{corollary}

The idea is that some level of a breadth-first search tree induces 
a graph of chromatic number at least $k + 1$. Such a graph 
contains an odd cycle of length at least $k(g - 2) + 1$ with 
a chord. We then apply the colouring lemma to deduce that 
$G$ contains cycles of $k(g - 2)$ consecutive lengths. In a sense, 
this generalizes a result of Gy\'{a}rf\'{a}s who showed that 
if $G$ has chromatic number at least $2k + 2$, then $G$ 
contains cycles of $k$ distinct odd lengths. These ideas may 
also be used to give a relatively short proof that 
$r(C_{2k+1},K_{n}) \ll n^{1 + 1/(k+1)}$. However, better results 
are easily obtain --- for example, it is possible to show that 
$r(C_{2k+1},K_{n}) \ll n^{1 + 1/(k+1)}(\log_{2}n)^{-1/k}$.

\bigskip

The cycles we have obtained are all very close together in the sense that they share many vertices. H\"{a}ggkvist and Scott~\cite{HS1} asked if it was possible, under an appropriate bound on the size of the graph, to find disjoint cycles of $k$ consecutive even lengths. This question also remains open, noting that a bound of order at least $k^{2}$ on the average degree would be required for disjoint cycles of $k$ consecutive even lengths. This is shown, for example, by $K_{l,n-l}$ where $l < 2k + k(k-1)/2$ and $n$ is sufficiently large. 

\medskip

We remark that from Theorem 1 we may obtain a result  on extremal numbers for even cycles,
 that slightly improves the result obtained by Bondy and Simonovits~\cite{BS} (see Corollary~9). From their paper, it follows that a graph of order $n$ and size at least $90kn^{1 + 1/k}$ contains a cycle of length $2k$. Two simple lemmas are required before proving our result. The first lemma is a special case of a lemma of Kostochka and Pyber~\cite{KP}. 
\begin{lemma}
 Let $G$ be a graph of order $n$ and size at least $cn^{1 + 1/k}$, where $c \geq 1$. Then $G$ contains a subgraph of average degree at least $c$ and radius at most $k$.
\end{lemma}
\begin{proof} We may assume that $G$ has minimum degree at least $cn^{1/k}$. Let $v_{0}$ be an arbitrary vertex in $G$ and define $H_{i}$ to be the subgraph of $G$ induced by vertices at distance at most $i$ from $v_{0}$. Define $r = \min\{i : e(H_{i}) \geq \frac{1}{2}c|H_{i}|\}$. Clearly $e(H_{i}) \geq \frac{1}{2} cn^{1/k} |H_{i-1}|$ for all $i$ and so, by definition of $r$, $|H_{r-1}| > n^{1/k}|H_{r-2}|$ which gives
$|H_{r-1}| > n^{(r-1)/k}$. Since $|H_{r-1}| \leq n$, $r - 1 < k$ and so $r \leq k$ and $H_{r}$ is the desired subgraph.
\end{proof} \medskip
\begin{lemma}
 Let $G$ be a graph of order $n$ with $e(G) \geq 2n^{1 + 1/k}$ where $k \geq 2$. Then $G$ has girth at most $2k+1$.
\end{lemma}
\begin{proof}
 By Lemma 6, $G$ has a subgraph $H$ of radius at most $k$ and average degree at least two. So $H$ contains a cycle and some cycle in $H$ has length at most $2k+1$.
\end{proof}

\medskip

In particular, if $G$ is bipartite and $k = 2$ in Lemma 7, then $G$ has a cycle of length four. For comparison, a standard result states that a graph which has at least $n^{3/2}/2 + n/4$ edges contains a cycle of length four.
We are able, using Theorem 1, to show the existence of longer even cycles:

\begin{theorem}
Let $G$ be a bipartite graph of order $n$ and girth $g$, and of size at least $4\lceil 2(k-1)/(g-2) \rceil n^{1 + 1/k}$, where $k \geq 2$ is an integer. Then $G$ has a cycle of length $2k$. 
\end{theorem}

\begin{proof}
By Lemma 7, either $G$ contains a cycle of length $2k$ or $g < 2k$ and $k \geq 3$. In the latter case, $\lceil 2(k-1)/(g-2) \rceil \geq 2$. 
Lemma 6 shows that $G$ contains a subgraph of average degree at least $4 \lceil 2(k-1)/(g-2) \rceil$ and of radius at most $k$. By Theorem 1, there are cycles of at least $k - 1$ consecutive even lengths in $H$, the shortest length being at most $2k$. So one of these cycles must have length exactly $2k$. 
\end{proof}

\medskip

 As a corollary to Theorem 8, we slightly improve the result of Bondy and Simonovits~\cite{BS}:

\begin{corollary}
Let $G$ be a graph of order $n$ and size at least $8(k-1)n^{1 + 1/k}$, where $k \geq 2$. Then $G$ contains a cycle of length $2k$.
\end{corollary}

\medskip

\noindent{\bf Acknowledgements}

\smallskip

I would like to thank Andrew Thomason for his many helpful suggestions.

\end{document}